\documentclass[oneside,english]{amsart}
\usepackage[T1]{fontenc}
\usepackage[latin9]{inputenc}
\usepackage{babel}
\usepackage{amsthm}
\usepackage{amstext}
\usepackage{amssymb}
\usepackage{esint}
\usepackage[unicode=true]
 {hyperref}
\usepackage{breakurl}

\makeatletter
\numberwithin{equation}{section}
\numberwithin{figure}{section}

\@ifundefined{definecolor}
 {\usepackage{color}}{}

\@ifundefined{definecolor}{\usepackage{color}}{}
\usepackage{amsfonts}\usepackage{babel}
\providecommand{\U}[1]{\protect\rule{.1in}{.1in}}

\usepackage{babel}

\makeatother

\begin{document}

\title[Kalton's work on differentials of complex interpolation]{An introduction to Nigel Kalton's work on differentials of complex
interpolation processes for K\"othe spaces}

\author{Michael Cwikel \and Mario Milman \and Richard Rochberg}

\address{Cwikel: Department of Mathematics, Technion - Israel Institute of
Technology, Haifa 32000, Israel }

\address{Milman: Department of Mathematics, Florida Atlantic University, Boca
Raton, FL 33431, USA }

\address{Rochberg: Department of Mathematics, Washington University, St. Louis,
MO, 63130, USA }

\email{mcwikel@math.technion.ac.il }

\email{mario.milman@gmail.com}

\email{rr@math.wustl.edu}

\thanks{The first author's work was supported by the Technion V.P.R.\ Fund
and by the Fund for Promotion of Research at the Technion. The second
author's work was partially supported by a grant from the Simons Foundation
(\#207929 to Mario Milman). The third author's work was supported
by the National Science Foundation under Grant No. 1001488.}

\maketitle
Our aim in this note is to offer one kind of introduction to Nigel
Kalton's remarkable paper \cite{K3}, and to share a few thoughts
about possible further sequels to it. We hope to at least capture
something of the spirit of the paper. Inevitably we will bypass many
of its subtleties. We will oversimplify and ignore details. 

At first sight, the main topics of \cite{K3} and the objects which
arise within them may seem quite exotic and even maybe ``far fetched''.
But it turns out that there are connections with and applications
to quite a range of other topics in analysis. We shall be able to
at least briefly hint at some of these below. Of course this is not
the only instance where Nigel Kalton's bold and deep explorations
along paths far from the ``beaten track'' have bounced back with
unexpected implications in more familiar settings. 

For other kinds of introductions to the same paper, we warmly encourage
the interested reader to consult the surveys \cite{Godf07} and \cite{GodfG2012}.
In \cite{Godf07} one can find a very extensive discussion of a wide
range of Nigel's research. In particular, its Section 5 deals, among
other things, with the material that we discuss here, as does \cite{GodfG2012}.
Each of these surveys provides many very illuminating insights and
offers quite different perspectives from ours, and mentions quite
a number of those connections with other topics to which we just alluded. 

Of course we also warmly recommend a somewhat longer and more detailed
account of these matters written by Nigel himself, together with Stephen
Montgomery-Smith in the later sections of their survey \cite{KM-S}. 

\medskip{}

\textbf{\textit{Acknowledgements:}} We thank Yoav Benyamini and David
Yost for some helpful correspondence.

\section{Interpolation}

\subsection{Classical Interpolation}

Classically, interpolation methods for Banach spaces are techniques
for starting with a pair of spaces, $(A_{0},A_{1})$, and constructing
an \textit{interpolation space} $A_{\ast}$ with the property that,
if a linear operator $T$ is bounded on $A_{i},i=0,1$ then $T$ is
also bounded on $A_{\ast}$.

One classical way of doing this is Alberto Calder\'on's \textit{complex
method of interpolation} (cf.~\cite{Ca}) which goes roughly as follows:
With ($A_{0},A_{1})$ given, define $\mathcal{F}=\mathcal{F}(A_{0},A_{1})$
to be the space of holomorphic vector valued functions, $F(z)$, defined
on the strip $S=\left\{ z\in\mathbb{C}:0\leq\operatorname{Re}z\leq1\right\} $
for which 
\[
\left\Vert F\right\Vert _{\mathcal{F}}=\max_{u=0,1}\sup\left\{ \left\Vert F(u+iv)\right\Vert _{A_{u}}:-\infty<v<\infty\right\} <\infty.
\]
 For $0<t<1$ one then defines the \textit{complex interpolation space
}$A_{t}$\textit{ }and its norm by 
\begin{align}
A_{t} & =\left[A_{0},A_{1}\right]_{t}=\left\{ F(t):F\in\mathcal{F}\right\} \nonumber \\
\left\Vert v\right\Vert _{A_{t}} & =\inf\left\{ \left\Vert F\right\Vert _{\mathcal{F}}:F\in\mathcal{F},\, F(t)=v\right\} ,\label{eq:Infimum}
\end{align}
 or, equivalently, $A_{t}$ is identified isometrically with the quotient
space 
\[
A_{t}=\mathcal{F}\diagup\left\{ F\in\mathcal{F}:F(t)=0\right\} .
\]
 The basic interpolation theorem is that if $T$ is bounded on $A_{i}$,
$i=0,1$ then $T$ is bounded on $A_{t}.$

\subsection{Differentials and Commutators}

Associated to the construction of $A_{t}$ is a special\textit{ ``differential''
map }$\Omega=\Omega(A_{0},A_{1},t)$ from $A_{t}$ into $A_{0}+A_{1}$.
For each $v\in A_{t},$ $\Omega(v)$ is the derivative at $t$ of
the function in $\mathcal{F}$ which attains the infimum in (\ref{eq:Infimum})
for defining $\left\Vert v\right\Vert _{A_{t}}.$ That is, define
$F_{t,v}$ by 
\[
F_{t,v}\in\mathcal{F}\text{, }v=F_{t,v}(t),\text{ }\left\Vert v\right\Vert _{A_{t}}=\left\Vert F_{t,v}\right\Vert _{\mathcal{F}},
\]
 and set 
\begin{equation}
\Omega v=F_{t,v}^{\prime}(t).\label{eq:GiveNumber}
\end{equation}
(A sample of the details we are ignoring here is the issue of whether
the above-mentioned infimum is attained and, if so, whether $F_{t,v}$
is unique, and what to do otherwise.) 

Operators such as $\Omega$ are the main topic of Nigel's paper. It
is easy to see that they are homogeneous, i.e., $\Omega(\alpha f)=\alpha f$
for all $\alpha\in\mathbb{C}$. It is not much harder to see that
for some constant $C$ and for all $f,g\in A_{t}$ 
\begin{equation}
\left\Vert \Omega(f+g)-\Omega(f)-\Omega(g)\right\Vert _{A_{t}}\leq C(\left\Vert f\right\Vert _{A_{t}}+\left\Vert g\right\Vert _{A_{t}}).\label{eq:RichardsNew2}
\end{equation}
Even though these operators are generally unbounded and nonlinear,
they interact well with the interpolation process. If $T$ is a linear
operator bounded on $A_{i},$ $i=0,1$ then, in addition to the boundedness
of $T$ on $A_{t},$ we have a commutator theorem (cf.~\cite{RW}):
there is a $C$ so that for all $v\in A_{t}$ 
\begin{equation}
\left\Vert \left[\Omega,T\right]v\right\Vert =\left\Vert \Omega(T(v))-T(\Omega(v))\right\Vert \leq C\left\Vert v\right\Vert .\label{eq:FirstCommutator}
\end{equation}

Here are some examples in a classical context of some other operators
which have similar properties to those of $\Omega$. Let $L^{2}=L^{2}(\mathbb{T},\mu)$
be the Lebesgue space on $\mathbb{T}=\left\{ e^{it}:0\le t<2\pi\right\} $
where $\mu$ is arc length measure. The role of the operator $T$
will be played by $P:L^{2}\to L^{2}$ which is the orthogonal projection
of $L^{2}$ onto the Hardy space, $H^{2},$ the subspace of $L^{2}$
consisting of functions whose Fourier coefficients with negative indices
vanish; $H^{2}=\left\{ f\in L^{2}:\left\langle f,e^{int}\right\rangle _{L^{2}}=0,\text{ }n<0\right\} .$ 

For each $f\in L^{2}$ and each $\tau=e^{it}\in\mathbb{T}$ set 
\begin{align*}
\left(\Omega_{1}f\right)(\tau) & =f(\tau)\log(1-\tau),\\
\left(\Omega_{2}f\right)(\tau) & =f(\tau)\log\frac{\left\vert f(\tau)\right\vert }{\left\Vert f\right\Vert _{L^{2}}},\text{ }\\
\left(\Omega_{3}f\right)(\tau) & =f(\tau)\log\mu\left(\left\{ \sigma\in\mathbb{T}:\left\vert f(\sigma)\right\vert >\left\vert f(\tau)\right\vert \right\} \right).
\end{align*}
(In fact the above formula for $\Omega_{3}f$ has to be replaced by
a more elaborate variant when $\left|f\right|$ assumes any constant
values on sets of positive measure.) Each of these operators has a
bounded commutator with $P$. That is, for each $i=1,2,3$ there is
a $C$ so that for all $f$ in $L^{2}$ 
\[
\left\Vert \lbrack\Omega_{i},P]f\right\Vert =\left\Vert (\Omega_{i}P-P\Omega_{i})f\right\Vert \leq C\left\Vert f\right\Vert .
\]
Note that none of the $\Omega_{j}$ are bounded on $L^{2}$. Also
$\Omega_{2}$ and $\Omega_{3}$ are nonlinear and yet we are claiming
linear space estimates. In fact $\Omega_{1}$ and $\Omega_{2}$ are
both differential maps $\Omega=\Omega(A_{0},A_{1},t)$ for suitable
choices of $A_{0}$, $A_{1}$ and $t$.

\subsection{More Generally}

There are other methods for constructing interpolation spaces; some
of the classical ones are discussed in, for example, \cite{BS}, \cite{BL},
\cite{BK}, and the earlier sections of \cite{KM-S}. Associated with
many of those methods are operators, similar to the special map $\Omega(A_{0},A_{1},t)$
defined above, but obtained by quite different constructions. For
example the operator $\Omega_{3}$ mentioned above, can be obtained
via an analogue of $\Omega(A_{0},A_{1},t)$ for one of the versions
of the method of real interpolation. $\Omega_{1}$, $\Omega_{2}$,
$\Omega_{3}$ and other such operators are discussed in \cite{CJMR},
\cite{JRW}, \cite{K1}, \cite{KM-S} (in its later sections) \cite{KM},
\cite{RW} and other places. 

Work to develop a more unified theory of such operators is in \cite{CCS},
\cite{C}, \cite{CKMR} and \cite{MR}. Some discussion of the applicability
of these constructions is in \cite{I} and \cite{R5}.

In the paper \cite{K3} Nigel pursues a different direction; he focuses
on complex interpolation and on a class of Banach spaces that are
amenable to more refined analysis. He can then explore more deeply
the relation between the interpolation construction and the associated
differentials. In fact, in\ \cite{K3} he shows that for Köthe function
spaces, one can develop a systematic theory for differentials $\Omega$,
and that the theory has interesting applications. We will discuss
that in the next section.

In recent years our understanding of interpolation theory has expanded
and several very interesting new interpolation constructions have
been introduced, for example in \cite{CS}, \cite{S1} and in many
papers by Zbigniew Slodkowski. (One might begin by looking at Slodkowski's
papers \cite{Slodkowski1988-A} and \cite{Slodkowski1988-B} and then,
in each case proceeding to the two subsequent similarly titled papers
which are their respective sequels.) These newer methods focus less
on boundedness results for linear operators and more on understanding
the role of convexity in Banach space theory; particularly the relation
with maximum principles and differential inequalities. This focus
on convexity and its variants (pseudoconvexity, quasiconvexity, quasi-affine
functions,...) in the theory of linear spaces was a major theme in
Nigel's research programs and his ideas in \cite{K3} resonate with
these newer views of interpolation. There is some brief discussion
of this in some of the comments which we offer in the final section.

\section{Kalton's Paper}

The work in \cite{K3} has strong connections with earlier results
in \cite{K1}. Its point of departure is a ``reasonable'' topological
space $S$ equipped with a $\sigma$-finite Borel measure $\mu$.
Nigel focuses on complex interpolation of a large class of Banach
spaces, K\"othe function spaces, whose elements are Borel measurable
functions on $S$. This class includes many Banach spaces which arise
naturally in various contexts and, consequently, the results of \cite{K3}
have quite a number of interesting applications. Within this class
there are two basic tools which are not available for general families
of Banach spaces. First, each space automatically carries with it
a rich natural family of bounded maps, namely multiplication by bounded
functions. Second, the underlying assumptions in this context ensure
that the dual of each space is also a space of functions on the same
underlying set and one can invoke the Lozanovsky's duality theorem
\cite{LozanovskiiG1969} (see also \cite{GillespieT1981} and \cite{ReisnerS1988})
and the associated factorization theory for functions.

\textbf{\textit{Remark:}}\textit{ Given that Lozanovsky's above-mentioned
result plays such a crucial and recurring role in what we are discussing
here, we see fit to mention the books \cite{Lozanovskii2012b,LozanovskiiG2000,LozanovskiiG2012a}
where the reader may discover other ideas of Lozanovsky. Some, maybe
even many of these may yet remain to be brought to fruition.}

\subsection{The Setup }

As we shall see, the main result of \cite{K3} follows from an elaborate
study of the interplay between several kinds of mappings and functionals,
in particular, \textit{derivations}, \textit{centralizers} and \textit{indicators},
which are defined on various Köthe function spaces, or other subsets
of the space $L^{0}$ of all measurable functions on the underlying
measure space.

To each space $A_{t}$ one can associate a new Banach space, Nigel's
\textit{derived space} $dA_{t}$ of couples $(u,v)$ in $A_{t}\times\left(A_{0}+A_{1}\right)$
for which the following norm is finite: 
\begin{equation}
\left\Vert (u,v)\right\Vert _{dA_{t}}=\inf\left\{ \left\Vert F\right\Vert :F\in\mathcal{F};\text{ }F(t)=u,F^{\prime}(t)=v\right\} .\label{eq:ThisNorm}
\end{equation}
It is relatively straightforward to see that this space coincides
with the \textit{twisted sum} (or \textit{twisted direct sum}) $A_{t}\oplus_{\Omega}A_{t}$
which is the space of pairs $(u,v)$ for which the following functional
\begin{equation}
\left\Vert (u,v)\right\Vert _{A_{t}\oplus_{\Omega}A_{t}}=\left\Vert u\right\Vert _{A_{t}}+\left\Vert v-\Omega u\right\Vert _{A_{t}}\label{eq:NewTwist}
\end{equation}
is finite. This functional is in fact a quasi-norm (in view of (\ref{eq:RichardsNew2}))
and it is equivalent to the norm (\ref{eq:ThisNorm}). (Cf.~the proof
of Lemma 2.9 of \cite[p.~325]{RW}. Note that here $\Omega$ is
the particular map $\Omega(A_{0},A_{1},t)$ defined by (\ref{eq:GiveNumber}).)
Using this fact it is not hard to check that the commutator estimate,
(\ref{eq:FirstCommutator}) for $T,$ is equivalent to knowing that
the map of $(u,v)$ to $(Tu,Tv)$ is bounded on $dA_{t}$ or, equivalently,
on $A_{t}\oplus_{\Omega}A_{t}.$ In fact Nigel also deals with these
notions in a broader context, via certain mappings $\Omega:A\to L^{0}$
which he calls \textit{derivations}, which generalize $\Omega(A_{0},A_{1},t)$.
(Twisted sums arise in still more general contexts and in some of
them, the quasi-norm (\ref{eq:NewTwist}) may fail to be equivalent
to a norm. See e.g.~\cite[p.~4]{GodfG2012} and, for further details,
Chapter 16 of \cite{BenLin} and the references therein.)

The basic question of \cite{K3} is: Does every such $\Omega$ arise
in this way? That is, given a Köthe function space $A$ and an $\Omega$
which satisfies various natural conditions, is there a couple $(A_{0},A_{1})$
and value of $t$ so that $A=A_{t}$ and $\Omega=\Omega(A_{0},A_{1},t)?$

To approach that question we first identify three natural necessary
conditions which $\Omega$ must satisfy. For any function $b$ let
$M_{b}$ be the operator of multiplication by $b.$

If $\Omega=\Omega(A_{0},A_{1},t)$ for some $t\in(0,1)$ and some
Köthe function spaces $A_{0}$, $A_{1}$ and $A=\left[A_{0},A_{1}\right]_{t}$
then it is not difficult to show that there must exist a positive
constant $\rho(\Omega)$ such that, for all $u,v\in A$, all $b\in L^{\infty}$,
and all $\alpha\in\mathbb{C}$, the map $\Omega$ satisfies 
\begin{equation}
\left.\begin{array}{ll}
(i) & \Omega(\alpha u)=\alpha\Omega(u)\\
(ii) & \Omega\left(\mathcal{B}_{A}\right)\mbox{ is bounded in }L^{0}\\
(iii) & \left\Vert \left[\Omega,M_{b}\right]u\right\Vert _{A}\leq\rho(\Omega)\left\Vert b\right\Vert _{\infty}\left\Vert u\right\Vert _{A}\,\,\,\,
\end{array}\right\} \label{eq:NewCentralizer}
\end{equation}
 (Here $L^{0}$ is equipped with its usual topology defined via convergence
in measure on sets of finite measure, and $\mathcal{B}_{A}$, as usual,
denotes the unit ball of $A$.) 

In particular, the third requirement in (\ref{eq:NewCentralizer})
arises because the $M_{b}$ are automatically bounded on all Köthe
function spaces, hence, for $\Omega$ to come from interpolation it
must satisfy commutator estimates with $M_{b}.$ 

Nigel uses the term \textit{centralizer} (or sometimes \textit{homogeneous
centralizer}) to describe any abstract map $\Omega:A\to L^{0}$ satisfying
(\ref{eq:NewCentralizer}) where $A$ is some K\"othe function space
(and there is no mention of any $A_{0}$ or $A_{1}$). He had already
considered similar maps in \cite{K1} and \cite{K2} (using the same
terminology but without imposing condition (ii)) and showed that they
are automatically also derivations in the more abstract sense alluded
to above.

\subsection{The Results}

The main result of \cite{K3} (Theorem 7.6 on page 511) is, roughly,
that if $\Omega$ satisfies (\ref{eq:NewCentralizer}) then it can
be obtained, to within a certain natural kind of equivalence, by complex
interpolation. That is, given $A$ and $\Omega$ one can select $A_{0}$,
$A_{1}$ and $t$ so that $A=A_{t}$ and $\Omega$ is equivalent to
$\Omega\left(A_{0},A_{1},t\right)$.

The proof involves an interesting associated construct, motivated
by Gillespie's alternative proof \cite{GillespieT1981} of the Lozanovsky
factorization theorem, namely the \textit{indicator} of the Köthe
function space $A$. This is the functional $\Phi_{A}$ defined initially
only on those positive elements $f$ of $L^{1}$ for which it is finite,
by 
\[
\Phi_{A}(f)=\sup_{x\in A.\left\Vert x\right\Vert _{A}\le1}\int_{S}f\text{ }\log\left\vert x\right\vert d\mu.
\]
(In other papers it is sometimes called the \textit{entropy function}.)
Nigel extends the proof of Lemma 3 of \cite{GillespieT1981} to obtain
that the supremum in the above formula is attained, for any given
$f$, by a positive function $x=x_{f}$ determined via the Lozanovsky
factorization theorem. This will enable him, at a rather later stage,
to extend the definition of $\Phi_{A}$ to a larger class of complex
valued functions $f$ by setting $\Phi_{A}(f)=\int_{S}f\log x_{\left|f\right|}d\mu$.
(This is ultimately done in Lemma 5.6 (on p.~499), but note that
there are misprints in the formula for $\Phi_{X}(f)$ on the third
line of the statement of that lemma, namely, the integral sign and
$"d\mu"$ have been omitted.) In parallel with his study of the indicator
functional $\Phi_{A}$, Nigel studies another more general class of
``indicator-like'' functionals $\Phi$ defined on suitable subsets
$\mathcal{I}$ of non-negative functions in $L^{1}$. His definition
of these functionals requires $\Phi(f)$ to be real when $f$ is real
valued, and also to satisfy certain continuity conditions. Furthermore,
using his notation 
\[
\Delta_{\Phi}(f,g)=\Phi(f)+\Phi(g)-\Phi(f+g)
\]
he requires that 

\begin{equation}
\Phi(\alpha f)=\alpha\Phi(f)\,\,\mbox{ \ensuremath{\forall}\ensuremath{\ensuremath{\alpha\ge}0}}\label{eq:indicatorONE}
\end{equation}
and also that, for some positive constant $\delta(\Phi)$ and all
$f,g\in\mathcal{I}$, 
\begin{equation}
0\le\Delta_{\Phi}\left(f,g\right)\le\delta(\Phi)\left(\left\Vert f\right\Vert _{L^{1}}+\left\Vert g\right\Vert _{L^{1}}\right).\label{eq:indicatorTWO}
\end{equation}

For each Köthe function space $A$, the indicator $\Phi_{A}$ has
all these properties, and, conversely, any $\Phi$ having all these
properties and also satisfying $\Delta_{\Phi}(f,g)\le\Delta_{\Phi_{L^{1}}}(f,g)$
for all $f,g\in\mathcal{I}$ is necessarily the indicator of some
Köthe function space $A$.

We will now become even more informal. Much of the technical work
in \cite{K3} is done using the indicators.

This \textquotedblleft{}change of variable\textquotedblright{}, to
working with the functionals $\Phi$ rather than the interpolation
spaces, \textquotedblleft{}linearizes\textquotedblright{} the interpolation
process: The indicator of $A_{t}=\left[A_{0},A_{1}\right]_{t}$ is
given by the formula:

\begin{equation}
\Phi_{A_{t}}=(1-t)\Phi_{A_{0}}+t\Phi_{A_{1}}\label{eq:NewLinear}
\end{equation}
and the Lozanovsky factorization for each K\"othe function space
$A$ can be expressed by the formula 
\[
\Phi_{L^{1}}=\Phi_{A}+\Phi_{A^{\ast}}.
\]

The fundamental technical conclusion of Nigel's paper is that one
can close the loop. Given a K\"othe function space $A$ and a $t$,
one can find $A_{0}$ and $A_{1}$ so that the indicator of $A$ satisfies
$\Phi_{A}=(1-t)\Phi_{A_{0}}+t\Phi_{A_{1}}$. It then follows that
$A$ is the interpolation space $A_{t}$. 

For the final steps towards his main result, Nigel has to reveal and
exploit a connection between centralizers and indicators. One should
keep in mind here that, while each K\"othe function space $A$ has
a unique indicator, there are infinitely many different centralizers
$\Omega$ which can be defined on $A$. Given any one such centralizer
$\Omega$, Nigel (as he already did in \cite{K1}) uses it to define
an auxiliary centralizer $\Omega^{\left[1\right]}$ defined on the
positive functions in $L^{1}$. Once more, ideas associated with the
Lozanovsky factorization play a central role, i.e., for each non negative
$x$ in $\mathcal{B}_{L^{1}}$,

\[
\Omega^{[1]}(x):=\Omega(a)a^{\ast}
\]
where $x=aa^{\ast}$ with $a$ in $A$, $a^{\ast}$ in the dual space
$A^{\ast}$, and with $\left\Vert a\right\Vert _{A}=\left\Vert a^{\ast}\right\Vert _{A^{*}}=1$.
Now it is possible to define a new functional $\Phi^{\Omega}$ on
a suitable subset of functions $f\in L^{1}$ by setting 
\[
\Phi^{\Omega}(f):=\int_{S}\Omega^{\left[1\right]}(f)d\mu.
\]
Results from \cite{K1} show that $\Phi^{\Omega}$ belongs to the
general class of ``indicator-like'' functionals mentioned above.
Therefore Nigel can apply his technical results about indicators:
Knowing that any indicator can be obtained by complex interpolation
insures that any centralizer can be obtained, to within an appropriate
equivalence, by complex interpolation.

\subsection{Furthermore}

Once the basic ideas are in place it is relatively straightforward
to obtain analogous results for rearrangement invariant spaces. That
is, if $A$ is a rearrangement invariant space and if given centralizers
or indicator functions interact appropriately with the linear operator
induced by rearrangements, then the new spaces constructed, $A_{0}$
and $A_{1},$ can be chosen to be rearrangement invariant. In particular,
by using rearrangement invariant spaces it is possible to obtain results
for the Schatten ideals. Those are spaces of compact operators on
Hilbert space which are normed by rearrangement invariant norms on
the operator's sequence of singular numbers. (Those are the numbers
which quantify the rate of approximation of a compact operator by
finite rank operators.)

Nigel also considers a converse question. If an operator is bounded
on a scale of spaces $A_{t},$ $0<t<1$ then the commutator with $\Omega,$
the associated derivation, is bounded on, say, $A_{1/2}.$ In the
other direction, what if we are given $A,$ $T,$ and a derivation
$\Omega,$ with both $T$ and $\left[\Omega.T\right]$ bounded on
$A_{1/2};$ must it follow that $T$ is bounded on the scale $\left\{ A_{t}\right\}$,
or, perhaps, at least for $t$ near $1/2?$ It is satisfying that
the answer is shown to be yes, at least in the case (in Theorems 9.7
and 9.8 \cite[pp.~524--525]{K3}) where $T$ is the (vector valued)
Riesz transform. It is intriguing that this answer can be used to
give an easy proof of a nontrivial result in harmonic analysis: A
fundamental fact about the ``good weights'' (i.e.~$A_{p}$ weights)
in the theory of singular integral operators is that their logarithms
are in the dual space of the Hardy space $H^{1}$, that is, they are
in $BMO;$ and, conversely, if the logarithm of a weight is in $BMO$
and has sufficiently small norm then it is a good weight. The first
of these facts can be derived using the theory of commutators that
we have been discussing; in fact the boundedness of $[P,\Omega_{1}]$
which we discussed earlier is a special case of that result. Nigel
can deduce the fact that exponentials of functions in $BMO$ are good
weights as a consequence (Corollary 9.9 \cite[p.~525]{K3}) of the
above-mentioned theorems). A classical presentation of these topics
and their uses is in Chapters 4 and 5 of \cite{SteinE1993}. Other
relations between the theory of commutators and classical analysis
are also developed in \cite{SteinE1993}. Other results indicating
the interplay between spaces which are analogues of $H^{1}$, and
of $BMO$ with classes of weights which are analogues of $A_{p}$,
can be found, for example, in \cite{BMR}. 

Theorems 9.7 and 9.8 also enable Nigel to obtain new results about
UMD-spaces.

\section{Looking Forward}

Some of the following observations seem to invite further research.
We will usually be brief and cryptic. Of course it is quite possible
that some of the questions that we ask here have been answered in
some publication which has not come to our attention. 
\begin{enumerate}
\item Some of the recent approaches to interpolation view interpolation
families $\left\{ A_{t}:0\leq t\leq1\right\} $ as ``geodesics''
in a space or ``manifold'' of possible Banach structures; \cite{CS},
\cite{CK}, \cite{R3}, \cite{S1} and \cite[Ch.~11]{S2}. From that
point of view, the construction of scales of complex interpolation
spaces is a method for constructing ``geodesics'' between two given
points; that is, solving a boundary value problem for geodesics. The
results in \cite{K3} show how to solve an initial value problem for
``geodesics''; that is, given a ``point'' (i.e., a K\"othe function
space $X$) and a ``tangent vector'' (i.e., a centralizer defined
on $X$), find a matching ``geodesic''. 
\item It is a theorem in classical analysis that, with $P$ denoting the
projection operator from our first example, then the bilinear map
$B(f,g)=f.Pg+g.Pf$, which at first glance maps $L^{2}\times L^{2}$
into $L^{1}$ actually has better properties; it maps into the real
variable Hardy space $\operatorname{Re}H^{1}.$ We mention two further
facts; some relations between them are developed in \cite{I} and perhaps
there is still more to learn. First, the properties of $B$ and related
maps can be used as the basis of a theory of compensated compactness
which is of great use when studying partial differential equations
(\cite{CLMS}). Second, this property of $B$ is a result in interpolation
theory. In particular it is equivalent, via a duality argument, to
commutator results of the type in (1) for $i=1.$ This point of view
is developed systematically by Nigel in \cite{K1} and \cite{K3}
where he develops a theory of a space $H_{\mathtt{sym}}^{1},$ \textit{the
symmetrized Hardy space, }a space which, in the classical case, is
closely related to the rearrangement invariant hull of the real variable
Hardy space. It is further shown, in \cite{K2}, that the Schatten
class associated to $H_{\mathtt{sym}}^{1}$ plays a fundamental role
in describing operator ideals (\cite{DFWW}). 
\item Interpolation constructions are generally both nonlinear and very
abstract. The passage to indicator functions replaces classical complex
interpolation with a explicit linear construction. It is not known
how indicators interact with other interpolation methods. 
\item Analytic semigroups of operators can be used to generate complex interpolation
families of Banach spaces (\cite[Sec.~4B]{RW}, \cite{R5}). The passage
from a scale of spaces to its differential seems analogous to passing
from a semigroup to its generator. Perhaps that analogy could be taken
further. There are interesting further thoughts developing this relation
in the final section of \cite{GodfG2012}. 
\item The theory of differentials and commutators associated to interpolation
extends to notions of higher differentials and associated commutators
(\cite{R4}, \cite{M}), but the formalism becomes relatively intricate.
Perhaps it is cleaner for Köthe function spaces. 
\item Also, some interpolation methods include an analysis of sub- and super-interpolation
families (\cite{R2}, \cite{S1}, \cite{CS}, \cite{CK}) which satisfy
various maximum or minimum principles. In some cases those are related
to curvature like expressions involving the second derivatives of
the norming function. Perhaps the higher differentials and commutators
provide a natural language in which to present such ideas. Perhaps
Köthe function spaces provide an area in which those ideas can be
explored more fully. 
\item Nigel's papers on commutators make very heavy use of the Lozanovsky
factorization. Are there implicit hints in his papers about how to
use the commutators to turn the machine around and run it in the other
direction? Could a rich theory of interpolation and commutators substitute
for Lozanovsky's duality and factorization theory in contexts far
removed from Köthe function spaces?
\item In a series of papers, H.~K\"onig and V.~Milman have studied operators
that satisfy certain functional equations. For example, they show
(cf.~\cite[Theorem 1]{KoMi}) that an operator $L:C^{1}(\mathbb{R})\rightarrow C(\mathbb{R})$
(not necessarily linear or continuous) that satisfies the classical
Leibniz rule for differentiating a product must be of the form $Lf=af^{\prime}+b\Omega f$,
where $a,b\in C(\mathbb{R})$, and $\Omega f=f\ln\left\vert f\right\vert $
. It could be of interest to study in detail the connection between
the K\"onig-Milman theory and the theory of commutators and its applications.
Here we make a few quick comments. K\"onig and Milman's characterization
allows them to conclude that for a Leibniz operator to act on higher
order spaces, e.g.~$L:C^{2}(\mathbb{R})\rightarrow C^{1}(\mathbb{R})$,
the ``cancellation condition'' $b=0$ must hold. Likewise, they
show that Leibniz operators on ``lower order spaces'', e.g.~$L:C(\mathbb{R})\rightarrow C(\mathbb{R})$,
must be of the form $Lf=c\Omega f$, for some $c\in C(\mathbb{R})$.
Moreover, the solution of functional equations of the form $Lf=g$,
has already appeared (as an auxiliary topic) in the study of higher
order logarithmic Sobolev inequalities (cf.~Feissner \cite[p.~58]{FeissnerG1978})
and commutator inequalities (cf.~\cite{CJMR}).
\item A key step in Nigel's journey to his main result is Theorem 6.6 on
p.~507. For a deeper understanding of this theory one might try to
determine whether the limiting case of either of these theorems holds,
i.e., when $\varepsilon=0$. Does the best value of the relevant
constant $C\left(\varepsilon\right)$ in this theorem necessarily
have to tend to $\infty$ as $\varepsilon$ tends to 0? The same question
could be asked perhaps more conveniently, with regard to Theorem 1.1
on p.~481, which is a finite dimensional \textquotedblleft{}model\textquotedblright{}
of Theorem 6.6 which Nigel formulates in his introduction to help
prepare the reader for what is to follow. In connection with this
question we also remark that the constant $\log2$ plays a special
role in these theorems. Furthermore, whenever $\Phi$ is the indicator
of a K\"othe space, the optimal (smallest) constant $\delta\left(\Phi\right)$
for which (\ref{eq:indicatorTWO}) holds satisfies $\delta\left(\Phi\right)\le\log2$.
Equality holds when the K\"othe space is $L^{1}$. For what other
spaces, if any, does equality hold?
\item As Nigel points out on p.~510, the centralizer $\Omega^{\left[1\right]}$
obtained when $\Omega=\Omega(A_{0},A_{1},t)$ is the same for all
values of $t\in(0,1)$. Can we somehow interpret this to mean that,
when $A_{0}$ and $A_{1}$ are both K\"othe function spaces, the
``geodesic curve'' $\left\{ [A_{0},A_{1}]_{t}:0\le t\le1\right\} $
which joins them is in some sense a ``straight line'' or has some
kind of zero ``curvature''. Within a whole ``manifold'' of Banach
spaces (which should all be compatibly contained in some Hausdorff
topological vector space) does it make sense to think of the K\"othe
function spaces, or Banach lattices as forming a ``flat'' ``submanifold''.
Is there some ``geometric'' way to characterize that ``submanifold''. 
\item Suppose that $A_{1}$, $A_{2}$, $A_{3}$ and $A_{4}$ are Banach
spaces which satisfy $A_{2}=\left[A_{1},A_{3}\right]_{\theta_{1}}$
and $A_{3}=\left[A_{2},A_{4}\right]_{\theta_{2}}$. Then Wolff's Theorem
\cite{WolffT1982} (see also \cite{JansonNilssonPeetre}) ensures
that $A_{2}=\left[A_{1},A_{4}\right]_{\alpha_{1}}$ and $A_{3}=\left[A_{1},A_{4}\right]_{\alpha_{2}}$
for suitable $\alpha_{1}$ and $\alpha_{2}$. In other words the ``geodesic
curves'' from $A_{1}$ to $A_{3}$ and from $A_{2}$ to $A_{4}$
can be ``glued together''| to form a single such ``curve'' from
$A_{1}$ to $A_{4}$. It is intriguing to note that, when all of the
above four spaces are K\"othe function spaces, this result emerges
from a trivial calculation using two instances of the formula (\ref{eq:NewLinear}). 
\item If $\Omega$ is a centralizer acting on some K\"othe function space,
then so is every constant multiple $r\Omega$ of $\Omega$, with the
constant $\rho(\Omega)$ replaced by $r\rho(\Omega)$. On the other
hand, if a centralizer is of the special kind $\Omega=\Omega(A_{0},A_{1},t)$
then it is not at all clear whether $r\Omega$ is exactly of this
kind, even if one renorms either or both of the spaces $A_{0}$, $A_{1}$.
Suppose that $A_{0}$ and $A_{1}$ are uniformly convex so that $\Omega(A_{0},A_{1},t)$
is uniquely defined, and that $B_{0}=A_{0}$ and $B_{1}=A_{1}$ but
with new norms $\left\Vert x\right\Vert _{B_{0}}=r_{0}\left\Vert x\right\Vert _{A_{0}}$
and $\left\Vert x\right\Vert _{B_{1}}=r_{1}\left\Vert x\right\Vert _{A_{1}}$
for each $x$ in $A_{0}$ or $A_{1}$ respectively. Then it follows
very easily that $\left[B_{0},B_{1}\right]_{t}=\left[A_{0},A_{1}\right]_{t}$
with $\left\Vert x\right\Vert _{\left[B_{0},B_{1}\right]_{t}}=r_{0}^{1-t}r_{1}^{t}\left\Vert x\right\Vert _{\left[A_{0},A_{1}\right]}$
for each $x\in\left[A_{0},A_{1}\right]_{t}$. But a simple calculation
gives an explicit formula which shows that $\Omega(A_{0},A_{1},t)$
cannot be a scalar multiple of $\Omega(B_{0},B_{1},t)$. However these
two maps are equivalent in the sense that the inequality 
\[
\left\Vert \Omega(A_{0},A_{1},t)x-c_{1}\Omega(B_{0},B_{1},t)x\right\Vert _{\left[A_{0},A_{1}\right]_{t}}\le c_{2}\left\Vert x\right\Vert _{\left[A_{0},A_{1}\right]_{t}}
\]
holds for suitable constants $c_{1}$ and $c_{2}$ and all $x\in\left[A_{0},A_{1}\right]_{t}$
These remarks show that it is not surprising that the precise formulation
of the main result Theorem 7.6 of \cite{K3} includes some requirements
on the size of the constant $\rho(\Omega)$. They also indicate that
there is no obvious way of obtaining a version of the theorem where
there is equality rather than merely equivalence of the associated
centralizers. 
\item In most of the natural applications of the results of \cite{K3},
the underlying measure space $\left(S,\mu\right)$ has a topological
structure, as is required by Nigel at the beginning of his exposition.
However, we suspect that most or maybe even all of the results of
\cite{K3} can be obtained in the context of a ``reasonable'' measure
space without any topology.
\item Since the paper \cite{K3} has such a wealth of ideas and powerful
methods, there are surely many more items that could be included here.
Should we have future thoughts about \cite{K3}, we may perhaps share
them with you, at least informally via the arXiv.\end{enumerate}

\end{document}